\documentclass[12pt]{amsart}
\usepackage{amsfonts, amssymb,amsmath, amsthm, amsxtra, latexsym, amscd}
\usepackage[latin1]{inputenc} 

\def\draft{\centerline{(Draft {\the \day}/{\the\month} \the \year.)}}

\theoremstyle{definition}
\newtheorem{theo+}    {Theorem}      [section]
\newtheorem{prop+}  [theo+]  {Proposition}
\newtheorem{coro+}  [theo+]  {Corollary}
\newtheorem{lemm+}  [theo+]  {Lemma}
\newtheorem{deep+}  [theo+]  {Deep Result}
\newtheorem{fact+}  [theo+]  {Fact}
\theoremstyle{definition}
\newtheorem{exam+}  [theo+]  {Example}
\newtheorem{rema+}  [theo+]  {Remark}
\newtheorem{defi+}  [theo+]  {Definition}
\newtheorem{xca+}[theo+]{Exercise}

\numberwithin{equation}{section}

\def\draft{\centerline{(Draft {\the \day}/{\the\month} \the \year.)}}
\bigskip

\def\refn#1.#2{\expandafter\def\csname#1\endcsname{[#2]}}
\def\refnr#1.{\csname#1\endcsname}

\def\Del{\Delta}

\def\a{\alpha}

\def\Claminv2{|C(\Lambda)|^{-2}}

\def\lam{\lambda}

\def\de{d\varepsilon}

\def\Aa2D{A^{\a,2}(D)}
\def\bAa2D{\overline{A^{\a,2}(D)}}
\def\Ab2D{A^{\beta,2}(D)}
\def\bAb2D{\overline{A^{\beta,2}(D)}}

\def\Norm#1_#2{\Vert#1\Vert_{#2}}

\def\phipl12{\phi_{p_{l_1}, p_{l_2}}}
\def\phip01{\phi_{p_{0}, p_{0}}}

\def\a{\alpha}

\def\Claminv2{|C(\Lambda)|^{-2}}

\def\Del{\Delta}

\def\sig{\sigma}

\def\lam{\lambda}

\def\tr{\operatorname{tr}}

\def\de{d\varepsilon}

\def\Aa2D{A^{\a,2}(D)}
\def\bAa2D{\overline{A^{\a,2}(D)}}
\def\Ab2D{A^{\beta,2}(D)}
\def\bAb2D{\overline{A^{\beta,2}(D)}}

\def\phipl12{\phi_{p_{l_1}, p_{l_2}}}
\def\phip01{\phi_{p_{0}, p_{0}}}

\def\bc{\mathbb C}
\def\br{\mathbb R}

\def\alg/{algebra} 
\def\Alg/{Algebra} 
\def\alt/{alternative} 
\def\anal/{analytic}
\def\analfunc/{\anal/\ \func/}
\def\Ans/{\it Answer. \normal}
\def\ass/{associative}
\def\nass/{non-\ass/}
\def\autom/{automorphism}
\def\homom/{homomorphism}
\def\isom/{isomorphism}
\def\bdd/{bounded}
\def\Bdd/{Bounded}
\def\bddsymdom/{bounded \sym/ \dom/}
\def\Cartdom/{Cartan \dom/}
\def\bdry/{boundary}
\def\bsd/{\bdd/ \symdom/}
\def\bv/{boundary value}
\def\cf/{{\it cf}\.}
\def\Cf/{{\it Cf}\.}
\def\charr/{character}
\def\coeff/{coefficient}
\def\comm/{commutative}
\def\cpct/{compact}
\def\compl/{complex}
\def\comp/{complex}
\def\Comp/{Complex}
\def\conf/{conformal}
\def\conj/{conjugate}
\def\conn/{connect}
\def\cont/{continuous}
\def\conv/{converge} 
\def\convc/{convergence}
\def\convt/{convergent}
\def\convx/{convex}
\def\coord/{coordinate}
\def\lcoord/{local coordinate}
\def\Corr/{Corresponding}
\def\corr/{corresponding}
\def\corrd/{correspond}
\def\cov/{covariant}
\def\decomp/{decomposition}
\def\deco/{decompose}
\def\diff/{different} 
\def\Diff/{Different} 
\def\dimn/{dimension} 
\def\distr/{distribution} 
\def\div/{diverge} 
\def\dom/{domain}
\def\eg/{\hbox{\it e.g}\.}
\def\eigenf/{eigen\-\func/}
\def\eigensp/{eigen\-space}
\def\eigenv/{eigen\-value}
\def\eq/{equation}
\def\equa/{equation}
\def\de/{\diff/ial \equa/}
\def\do/{\diff/ial operator}
\def\ode/{ordinary \de/}
\def\pde/{partial \de/}
\def\pdo/{partial \diff/ial operator}
\def\psdo/{pseudo \diff/ial operator}
\def\fin/{finite}
\def\Ex/{\it Example.\ \normal}
\def\Exnr#1/{\it Example #1.\ \normal}
\def\foll/{follow}
\def\follg/{following}
\def\Follg/{Following}
\def\func/{function}
\def\Func/{Function}
\def\Fonc/{Fonc\-tion}
\def\fonc/{fonc\-tion}
\def\Funk/{Funk\-tion}
\def\funk/{Funk\-tion}
\def\gen/{general}
\def\har/{harmonic}
\def\Hint/{\it Hint. \normal}
\def\hist/{historic}
\def\histcl/{historical}
\def\hol/{holo\-morphic}
\def\homog/{ho\-mo\-ge\-ne\-ous}
\def\hyp/{hyper\-bolic}
\def\hyperg/{hyper\-geometric}
\def\ie/{\hbox{\it i.e.}}
\def\iff/{if and only if}
\def\ineq/{inequality}
\def\infra/{{\it inf\-ra}}
\def\ultra/{{\it ult\-ra}}
\def\Inpart/{In particular}
\def\inpart/{in particular}
\def\instof/{instead of}
\def\interps/{interpolation space}
\def\interp/{interpolation}
\def\Interp/{Interpolation}
\def\interpr/{Interpretation}
\def\Intr/{Introduction}
\def\intv/{interval}
\def\inv/{invariant}
\def\invc/{invariance}
\def\Iowords/{In other words}
\def\iowords/{in other words}
\def\ipr/{inner product}
\def\irred/{irreducible}
\def\lb/{line bundle}
\def\lin/{linear}
\def\lhs/{left hand side}
\def\rhs/{right hand side}
\def\loc/{local}
\def\math/{mathematic} 
\def\mathcn/{\math/ian}
\def\manif/{manifold}
\def\meas/{measure}
\def\measl/{measurable}
\def\mero/{mero\-morphic}
\def\mon/{monomial}
\def\monog/{monogenic}
\def\mult/{multiple}
\def\multy/{multiply}
\def\multn/{multiplication}
\def\nas/{necessary and sufficient}
\def\nbd/{neighborhood}
\def\neg/{negative}
\def\nondeg/{nondegenerate}
\def\Oohand/{On the other hand}
\def\oohand/{on the other hand}
\def\Oonhand/{On the one hand}
\def\oonhand/{on the one hand}
\def\oper/{operator}
\def\orth/{ortho\-gonal}
\def\orthon/{ortho\-normal}
\def\otoh/{on the other hand}
\def\quat/{quaternion}
\def\pp/{\hbox{a. e.}}
\def\psh/{plurisubharmonic}
\def\pol/{polynomial}
\def\pot/{potential}
\def\pos/{positive}
\def\princ/{principle}
\def\prob/{probability}
\def\proj/{projective}
\def\projn/{projection}
\def\Proof/{\it Proof:\normal}
\def\Rem/{\it Remark\normal}
\def\Remnr#1/{\it Remark\ \normal #1. }
\def\rep/{representation}
\def\meta/{metaplectic representation}
\def\repr/{reproducing}
\def\reprker/{reproducing kernel}
\def\resp/{respective} 
\def\resply/{respectively}
\def\restr/{restriction}
\def\sa/{self-adjoint}
\def\st/{such that}

\def\sol/{solution}
\def\ru/{space}
\def\sph/{spherical}
\def\ssp/{sub\ru/}
\def\sym/{symmetric}
\def\Sym/{Symmetric}
\def\symb/{symbol}
\def\symbc/{symbolic}
\def\symdom/{\sym/ domain}
\def\symp/{symplectic}
\def\Theor#1/{\fet Theorem #1.\ \normal}
\def\Lem#1/{\fet Lemma #1.\ \normal}
\def\Lemma/{\fet Lemma.\ \normal}
\def\topl/{topology}
\def\topll/{topological}
\def\transf/{transform}
\def\transl/{translation}
\def\transfn/{transformation}
\def\transv/{transvectant}
\def\trig/{trigonometric}
\def\tril/{trilinear}
\def\trilf/{trilinear form}
\def\uhp/{upper halfplane}
\def\uhs/{upper halfspace}
\def\vb/{vector bundle}
\def\vf/{vector field}
\def\vsp/{vector space}
\def\wrt/{with respect to}
\def\Wlog/{Without loss of generality}
\def\a{\alpha}

\def\lam{\lambda}

\def\sig{\sigma}

\def\Ab/{Abel}
\def\Ban/{Banach}
\def\Bansp/{\Ban/ space}
\def\Belt/{Bel\-tra\-mi}
\def\Berg/{Berg\-man}
\def\Bern/{Ber\-nou\-lli}
\def\Berz/{Berezin}
\def\Bess/{Bessel}
\def\Cart/{Car\-tan}
\def\Cay/{Cay\-ley}
\def\CG/{Clebsch-Gordan}
\def\Cl/{Clifford}
\def\CR/{Cauchy-Rie\-mann}
\def\Dir/{Dirichlet}
\def\Eucl/{Euclide}
\def\F/{Fourier}
\def\Hank/{Hankel}
\def\Hankf/{\Hank/ form}
\def\Herm/{Hermite}
\def\Hilb/{Hilbert}
\def\Hilbs/{Hilbert space}
\def\Hilbsp/{Hilbert space}
\def\HS/{Hilbert-Schmidt}
\def\Lag/{La\-grange}
\def\Lap/{La\-place}
\def\LapBelt/{\Lap/-\Belt/}
\def\Leb/{Lebesgue}
\def\Marc/{Mar\-cin\-kie\-wicz}
\def\Moeb/{Moebius}
\def\Moebt/{Moebius transformation}
\def\Moebtransfn/{Moebius transformation}
\def\Pla/{Plan\-che\-rel}
\def\Poin/{Poin\-car\'e}
\def\Riem/{Rie\-mann}
\def\Riemn/{\Riem/ian}
\def\psRiemn/{pseudo-\Riem/ian}
\def\Riems/{Rie\-mann surface}
\def\Schroe/{Schr\"odinger}
\def\Weier/{Weier\-strass}
\def\anal/{analytic}
\def\bsd/{bounded symmetric domain  }
\def\bdd/{bounded}
\def\calc/{calculation}\def\conj{conjugate}
\def\calci/{calculating}\def\eg{e.g.}
\def\conj/{conjugate}
\def\deco/{decomposition}
\def\eg/{e.g.}
\def\fct/{function}
\def\gp/{group}
\def\hw/{highest weight}
\def\hwv/{highest weight vector}
\def\hwvs/{highest weight vectors}
\def\lw/{lowest weight}
\def\lwv/{lowest weight vector}
\def\lwvs/{lowest weight vectors}
\def\hds/{holomorphic discrete series}
\def\iff/{if and only if}
\def\inv/{invariant}
\def\irrde/{irreducible decomposition}
\def\meas/{measure}
\def\transf/{transform}
\def\rep/{representation}
\def\resp/{respectively}
\def\inters/{intertwines}
\def\interg/{intertwining}
\def\meta/{metaplectic representation}
\def\qu/{quaternion}
\def\rep/{representation}
\def\symdom/{ symmetric domain}
\def\st/{such that}
\def\shd/{subhead}
\def\transf/{transform}
\def\wrt/{with respect to}

\def\Norm#1#2#3{\Vert#1\Vert^{#3}_{{#2}¥}}

\def\tr{\operatorname{tr}}

\begin{document}

\def\bbK{\mathbb K}
\def\bKk{\mathbb K^k}
\def\bKn{\mathbb K^n}
\def\Cos^2{\text{Cos}^2}

\def\mP{\mathcal P}
\def\mF{\mathcal F}
\def\SdN{\mathcal S^N}
\def\eSdN{\mathcal S_N^\prime}
\def\Gc{G_{\mathbb C}}
\def\rqs{G\backslash K_{\mathbb C}}

\baselineskip 1.25pc


\def\psdo{\text{$\Psi D O$}}
\def\trr{\tr_{\text{res}}}
\def\tro{\tr_{\omega}}
\newcommand\dbar{\overline\partial}
\renewcommand\Im{\operatorname{Im}}
\hfuzz6pt

\title[Toeplitz and Hankel operators and Dixmier traces]
{Toeplitz and Hankel operators and Dixmier traces
on the unit ball of $\mathbb C^n$}

\author{Miroslav Engli\v s, Kunyu Guo and Genkai Zhang}

\address{Mathematics Institute~AS~\v CR, \v Zitn\'a~25, 11567~Prague~1, Czech
Republic {\rm and} Mathematics Institute, Silesian University, Na~Rybn\'\i\v
cku~1, 74601~Opava, Czech Republic}
\email{englis@math.cas.cz}
\address{Department of Mathematics, Fudan University, Shanghai 200433, P. R. China}
\email{kyguo@fudan.edu.cn}
\address{Department of Mathematics, Chalmers University of Technology and
G\"oteborg University, G\"oteborg, Sweden}
\email{genkai@math.chalmers.se}
\thanks{Research of M. Engli\v s supported by GA~\v CR grant
no.~201/06/128 and AV~\v CR Institutional Research Plan no.~AV0Z10190503,  K. Guo by NSFC(10525106) and  NKBRPC(2006CB805905) and 
G. Zhang by the Swedish Science Council (VR)
and SIDA-Swedish Research Links}

\keywords{Schatten - von Neumann classes, Macaev classes, trace, Dixmier trace,
Toeplitz operators, Hankel operators, pseudo-Toeplitz operators,
pseudo-differential operators, boundary CR operators, invariant Banach spaces}
\begin{abstract}
We~compute the Dixmier trace of pseudo-Toeplitz
operators on the Fock space. As~an application we find
a formula for the Dixmier trace of the product
of commutators of Toeplitz operators on the Hardy and weighted
Bergman spaces on the unit ball of~$\mathbb C^d$. This 
generalizes an earlier work of Helton-Howe
for the usual trace of the anti-symmetrization
of Toeplitz operators.
\end{abstract}

\maketitle

\section{Introduction}
In the present paper we will study the Dixmier
trace of a class of Toeplitz and Hankel operators
on the Hardy and weighted Bergman spaces on 
the unit ball of $\mathbb C^d$. We give
a brief account of our problem and explain 
some  motivations. 
Consider the Bergman space
$L^2_a(D)$ of holomorphic functions on
the unit disk $D$ in the complex plane. For a bounded function $f$
let $T_f$ be the Toeplitz operator on $L^2_a(D)$. It
is  well-known that for a holomorphic function $f$
the commutator $[T_f^\ast, T_f]$ is of trace
class and the trace is given by the square of the Dirichlet norm of~$f$, 
$$
\tr [T_f^\ast, T_f]=\int_D |f'(z)|^2 \, dm(z),
$$
which is one of the best known M\"obius invariant
integrals. This formula actually holds
for  Toeplitz operators on any Bergman
space on a bounded domain with the area
measure replaced any reasonable measure \cite{AFJP-crelle}.
There is a significant difference
between Toeplitz operators
on the unit disk and on the unit ball $B=B^d$ in
$\mathbb C^d$, $d>1$. Let $\mathcal L^p$  be the Schatten - von Neumann
class of $p$-summable operators. The commutator $[T_f^\ast, T_f]$ on the 
weighted Bergman space, say for holomorphic functions $f$
in a neighborhood of the closed the unit disk, 
is in the Schatten - von Neumann class $\mathcal L^{p}$,  for $p>\frac 12$
and is zero if it is in  $\mathcal L^{p}$,  for $p\le \frac 12$, $\frac 12$
being called the cut-off; on the Hardy space
$[T_f^\ast, T_f]$ can be in any 
 Schatten - von Neumann class $\mathcal L^{p}$,  for $p>0$;
see \cite{Peller-book} 
and  \cite{Rochberg-iumj82} for the case of Hardy space and 
 \cite{afp-ajm88} for the case of weighted Bergman space.
However for $d>1$, it is in $\mathcal L^{p}$ for $p >d$,
with $p=d$ being the cut off,
both on the weighted Bergman spaces and on the Hardy space.
Thus no trace formula was expected
for the commutators. Nevertheless Helton and Howe \cite{Helton-Howe-acta} 
were able to find an analogue of the previous formula. They
showed, for smooth functions $f_1, \cdots, f_{2d}$
on the closed unit ball, that
the anti-symmetrization $[T_{f_1}, T_{f_2},  \cdots T_{f_{2d}}]$ of 
the $2d$ operators $T_{f_1}, T_{f_2}, \cdots T_{f_{2d}}$ 
is of trace class and  found that
$$
\tr
[T_{f_1}, T_{f_2},  \cdots T_{f_{2d}}] =\int_B df_1\wedge df_2\cdots \wedge
df_{2d}. 
$$

On the other hand,  we observe that $[T_{f}, T_{g}]
$ is, for
smooth functions $f$ and $g$,
in the Macaev class $\mathcal L^{d, \infty}$ (which is an analogue
of the Lorentz space $L^{d, \infty}$), thus
the product of $d$ such commutators
$[T_{f_1}, T_{g_1}][T_{f_2}, T_{g_2}]\cdots
[T_{f_d}, T_{g_d}]$ is in  $\mathcal L^{1,\infty}$
and hence has a Dixmier trace. One~of the goals of the
present paper is to prove the following
formula for the Dixmier trace of this product of commutators:
$$
\tro[T_{f_1}, T_{g_1}]\cdots
[T_{f_d}, T_{g_d}]
=\int_S \{f_1, g_1\}
\cdots \{f_d, g_d\}.
$$
Here $ \{f, g\}$ is the Poisson bracket of $f$ and~$g$; its
restriction to the boundary $S$ of $B$ depends only on the boundary
values of $f$ and $g$ and can be expressed in terms
of the boundary $CR$ operators. This can be viewed as a generalization
of the Helton-Howe theorem. We apply our result also
to Hankel operators and obtain a formula for the Dixmier trace of the $d$-th
power of the square modulus of the Hankel operators $H_{\overline f}^\ast
H_{\overline f}$ for holomorphic functions~$f$. Namely
we have
$$
\tro |H_{\bar f}|^{2d}
= \tro ([T_{\bar f}, T_f]^{d})
= \int_{S}(|\nabla f|^2-|Rf|^2)^d .
$$
This provides
a boundary $\mathcal L^{d, \infty}$
result for the  Schatten-von~Neumann $\mathcal L^{p}$ ($p>d$)
properties of the square modulus of the Hankel operators 
(see \cite{AFJP-jlms}, \cite{Roch-Feld} 
and \cite{Zhu-ajm91}).  We mention also, besides
the above results on exact norms, there
are exact formulas proved by
Janson, Upmeier and Wallsten \cite{JUW-jfa}
on the Schatten - von Neumann $\mathcal L^p$-
norm of the Hankel operators on the unit circle for $p=2, 4, 6$,
and by Peetre \cite{Peetre-s4} on $\mathcal L^4$-norms
of Hankel forms on Fock spaces \cite{JPR-revista}. 

There has been an intensive study of Dixmier trace and residue
trace of pseudo-differential
operators, mostly on compact manifolds where the analysis
is relatively easier, see e.g. \cite{jfa-dixmier-psdo}
and \cite{ponge-arxiv-jul06} and references therein, thus
the Toeplitz  operators
on Hardy spaces on the boundary of a bounded strictly
pseudo-convex domains can be treated using
the techniques developed there. The Hankel and Toeplitz
operators on  Bergman spaces, generally speaking, behaves
rather differently from those on Hardy space, and 
the result of Howe \cite{Howe-jfa80} roughly speaking
proves that Toeplitz operators of certain classes
 can be treated similarly to those in Hardy space case
(also called the de Monvel - Howe compactification  \cite{Guillemin-ieot}).
Our result can thus be viewed a generalization
of the  compactification to weighted Bergman spaces
and an application of the  \cite{Connes-cmp88}
ideas of computing
Dixmier traces.
In  particular
our Theorem 4.1 are closely related
to the  results in \cite{jfa-dixmier-psdo} where the residue
trace of pseudo-differential operators
of certain class is computed; here we use
the Weyl transforms and they differ from 
pseudo-differential operators of lower order, so
that Theorem 4.1 can also be obtained from \cite{jfa-dixmier-psdo}
provided
one proves the the lower order terms are of trace
class.

In another paper we will study the Dixmier trace for Toeplitz operators on a
general strongly pseudo-convex domain. One~of the authors, G.~Zhang, would also
like to thank Professor Richard Rochberg and Professor Harald Upmeier for
introducing him   the work of Connes \cite[Chapter IV.2]{Connes-ncg} on Dixmier
traces of pseudo-differential operators.

\section{Toeplitz operators  on Bergman spaces and their
realization as pseudo-Toeplitz operators on Fock spaces}

Let $dm(z)$ be the Lebesgue measure on $\bc^d$ and
consider the weighted measure
$$
d\mu_{\nu}=C_{\nu}(1-|z|^2)^{\nu-d-1}dm(z),
$$
where $C_{\nu}$ is the normalizing constant to make $d\mu_\nu$ a probability
measure and \hbox{$\nu>d$}.
We~let $\mathcal H_{\nu}$ be the corresponding Bergman 
space of holomorphic functions on $B$. We will also
consider the Hardy space of
square integrable functions on $S$ which are
holomorphic on $B$. This can be viewed
as the analytic continuation of $\mathcal H_\nu$ at $\nu=d$.
Thus we assume throughout this paper that $\nu \ge d$.

Let $f$ be a bounded smooth function on $\bar B$, the closure
of $B$. The Toeplitz operator $T_f$ on $\mathcal H_\nu$
with symbol $f$ is defined by
$$
T_f  g= P (fg)
$$
where $P$ is the Bergman or the Hardy projection
for $\nu>d$ and $\nu=d$, respectively.

As~was shown by Howe \cite{Howe-jfa80} there is
a more flexible and effective
way of studying the spectral properties
of Toeplitz operators with smooth symbol,
by using the theories  of representations
of the Heisenberg group and  of pseudo-differential operators. 
We will adopt that approach. We will be very
 brief and refer to
 \cite{Howe-jfa80} and \cite[Chapter XII]{Stein-book2}
for details.
 So let $H_{n}=\mathbb C^d \times T$ 
be the Heisenberg group as in loc. cit.. The Heisenberg
group has an irreducible representation, $\rho$,  on the
Fock space $\mathcal F$ consisting
of entire functions $f$ on
$\mathbb C^d$ such
that
$$
\int_{\bc^d} |f(z)|^2 e^{-\pi |z|^2} dm(z) <\infty.
$$
The action of the Heisenberg group is explicitly
given as follows. For $w\in \mathbb C^d$ viewed 
as an element in $H_d$,
$$
\rho(w)f(w^\prime)= e^{-\pi/2 |w|^2 +\pi w^\prime \cdot\overline w}
 f(w^\prime -w),
$$
where $w^\prime\cdot\overline w$ is the Hermitian inner product on
$\mathbb C^d$. The action of $T$ is by rotation. 

Identifying the Lie algebra $\mathfrak h$
of the Heisenberg group with $\mathbb R^{2n}\oplus \mathbb R$
and thus $\mathbb R^{2n}$ with a subspace of the Lie algebra
we get an action of $\mathbb R^{2n}$ as holomorphic differential
operators on $\mathcal F$, which extends from $\mathfrak h$ to the whole
enveloping algebra $\mathfrak U(\mathfrak h)$ and
which will also be denoted by $\rho$.
In~particular, taking the basis elements $\partial_j=\partial/\partial w_j$ 
and $\dbar_j=\partial/\partial\overline w_j$ of $\mathbb R^{2n}$ we have
\begin{equation}
  \label{eq:ossci-1}
\rho(\partial_j) f(w) = -\partial_j f(w), \qquad
\rho(\dbar_j) f(w) = \pi w_j f(w).
\end{equation}
Let, following the notation in \cite{Howe-jfa80}, 
$\Delta\in\mathfrak U(\mathfrak h)$ be the element
$$
\Delta = \frac 12 (\partial_j\cdot\dbar_j+\dbar_j\cdot\partial_j).
$$
Then $\rho(\Delta)$ acts on $\mathcal F$ as a diagonal self-adjoint  operator
\cite{Howe-jfa80},
under the orthogonal basis $\{w^\alpha, \alpha=(\alpha_1, \cdots, \alpha_d)\}$,
viz
\begin{equation}
  \label{eq:ossci}
\rho(\Delta) w^{\alpha}   =-\pi\Big(|\alpha| +\frac d2\Big) w^{\alpha}.  
\end{equation}

Let $F(z)$ be a function on $\mathbb C^d$ (viewed
as a  function on the Heisenberg  group).
The Weyl transform $\rho(F)$ of $F$ is defined
by
$$
\rho(F)=\int_{\mathbb C^d} F(w) \rho(w) dm(w).
$$

To understand the operator theoretic properties of 
$\rho(F)$
we will need the Fourier transform of $F$. 
Let $\hat F$ be the (symplectic-) Fourier transform of
$F$
$$
\hat F(w^\prime)
=2^{-d} \int_{\bc^d} F(w) e^{\pi i \Im w^\prime\cdot\overline w } dm(w),
$$
and $F\ast G$ the  symplectic convolution
$$
F\ast G (w)=\int_{\bc^d} F(z) G(w-z) e^{\pi i \Im w\cdot\overline z} dm(z).
$$
We recall that
$$
\widehat {F \ast G}=  F \ast \hat  G
$$
and
$$ \rho(F) \rho(G) = \rho(F*G)  $$
for appropriate class of functions.
A well-known theorem of Calder\'on-Vaillancourt states
that if $\hat F$ and all its derivatives are bounded then
$\rho(F)$ can be defined as a bounded operator on $\mathcal F$.

We will need a finer class of symbols introduced by Howe, 
corresponding to the so-called pseudo-Toeplitz operators.
Let
$$\mathcal {PT}(m, \mu)
=\{F\in \mathcal S^\ast (\bc^d): 
|{\partial}^{\alpha}{\dbar}^{\beta}\hat F| 
\le C_{\alpha \beta}(1+|w|)^{m-\mu(|\alpha| +|\beta|)}
\}$$
and
$$\mathcal {PT}_{rad}(m, \mu)=\{F\in \mathcal {PT}(m, \mu):
\hat F=(1-g(|w|))\psi(\frac{w}{|w|}) |w|^m + D_1, \, D_1\in 
\mathcal {PT}(m-\mu, \mu) \}.
$$
Here $g$ is a smooth function on $\br$ such that
$0\le g(t)\le 1$ on $\mathbb R$, $g(t)=0$ for $|t|\ge 2$ and
$g(t) =1$ for $0\le t\le 1$. 

 For $F\in \mathcal {PT}_{rad}(m, \mu)$ we will call
\begin{equation}
  \label{eq:pr-sym}
\sig_{m}(F):=\psi(\frac{w}{|w|}) |w|^m  
\end{equation}
its principal symbol. It can be obtained, up~to the factor $|w|^m$,~by
$$
\psi(w)=\lim_{t\to \infty} t^{-m} \hat F (tw), \quad  w\in S.
$$
Following Howe~\cite{Howe-jfa80}, we~will call $\rho(F)$,
$F\in\mathcal{PT}(m,\mu)$, a~pseudo-Toeplitz operator of order $m$ and
smoothness~$\mu$. One~has \cite[Lemma~4.2.2]{Howe-jfa80}
\begin{equation}
 F\in\mathcal{PT}(m_1,\mu), \; G\in\mathcal{PT}(m_2,\mu)\;
 \implies F\ast G\in\mathcal{PT}(m_1+m_2,\mu).
\label{eq:pridana} \end{equation}

We will realize  the Toeplitz operators $T_f$
on $\mathcal H_{\nu}$ for $f$ on $B$
(or on $S$ for the Hardy space)
as Weyl transforms $\rho(F)$ of certain symbols $F$
on $\mathbb C^d$.  First we notice that
$$
e_{\beta}:=\Big(\frac{(\nu)_{|\beta|} }  {\beta!} \Big)^{\frac 12} z^{\beta}
$$
form an orthonormal basis of
$\mathcal H_{\nu}$, and so~do
$$
E_{\beta}:=\Big(\frac{1}{\pi^{|\beta|}\beta!}\Big)^{\frac 12} w^{\beta}
$$
for $\mathcal F$. (Here $(\nu)_j:=\nu(\nu+1)\dots(\nu+j-1)$
is the usual Pochhammer symbol.) Thus the map
\begin{equation}
  \label{eq:uni-ber-foc}
U: e_{\beta} \to E_{\beta}
\end{equation}
is an unitary operator.
  First we will find
the action of the elementary 
Toeplitz 
operators $T_{z^{\alpha}}$ under the intertwining map $U$.
\begin{lemm+} The operator $U T_{z^{\alpha}} U^\ast$ on $\mathcal F$ 
is given~by
\begin{equation}
   \label{eq:int-ta}
U T_{z^{\alpha}} U^\ast =
 {\rho(z)}^{\alpha}
 \rho\left( \pi^{|\alpha|}
 \Big(\nu-\frac d2 -\frac 1\pi \Del\Big)_{|\alpha|}  \right)^{-1/2}   
 \end{equation}
\end{lemm+}
This can be proved by direct computation. Indeed we have
$$
T_{z^{\alpha}} e_{\beta}
= \Big(\frac{(\beta)_{\alpha}}
{(\nu +|\beta|)_{|\alpha|}}\Big)^{\frac 12}
e_{\beta +\alpha},
$$
and the right hand side (\ref{eq:int-ta})
can be  easily computed by (\ref{eq:ossci-1}) and  (\ref{eq:ossci}).

By using the previous Lemma we have then the following result which was 
proved by Howe \cite[Proposition 4.2.3]{Howe-jfa80} in the case when $\nu=d+1$;
the general case of $\nu\ge d$ is essentially the same.

\begin{prop+} Let $f\in C^\infty(S)$ and let $\tilde f$
be a $C^\infty$ extension to $B$ and $T_{\tilde f}$ the Toeplitz
operator on $\mathcal H_\nu$. Then under the unitary equivalence 
of $\mathcal H_\nu$ and the Fock space $\mathcal F$ on $\mathbb C^d$, 
the~Toeplitz operators are pseudo-Toeplitz operators with radial asymptotic
limits $\mathcal {PT}_{rad}(0, 1)$.
More precisely, there exists $F\in\mathcal{PT}_{rad}(0,1)$ such that 
$UT_{\tilde f}U^*=\rho(F)$, and $f(\zeta)=\lim_{t\to\infty} \hat F(t\zeta)$
for each $\zeta\in S$.
\end{prop+}

\section{Schatten-von Neumann properties of pseudo-Toeplitz operators}
Recall that the Schatten-von Neumann class $\mathcal L^p$, $p\ge 1$,
consists of compact operators $T$ such that the eigenvalues $\{\mu_n\}$ 
of $|T|=(T^\ast T)^{\frac 12}$ are $p$-summable,
$\sum \mu_n^p <\infty$. In~particular
$\mathcal L^2$ is the Hilbert-Schmidt class,
$\mathcal L^1$  the trace class and $\mathcal L^\infty$
are the compact operators. For $1< p<\infty$, $1\le q\le \infty$, 
the~Macaev class $\mathcal L^{p, q}$ is obtained by the real interpolation
between $\mathcal L^1$ and $\mathcal L^\infty$.
However, we~will need the Macaev class $\mathcal L^{1,\infty}$, 
which consists \footnote{
Sometimes this ideal is denoted by $\mathcal L_\Omega$,
       the notation $\mathcal L^{1,\infty}$ being reserved for
       the (smaller) class of operators for which $\mu_n=O(1/n)$.
       Our notation follows Coones' book [9].
}
of all compact operators such that, if~$\mu_1\ge\mu_2\ge\dots$,
$$
\sum_{n=1}^N\mu_n =O(\log N).
$$

There exists a linear functional on the space
$\mathcal L^{1, \infty}$ that resembles the usual trace, called the
Dixmier trace. Its definition is rather
involved and we refer to \cite[Chapter~IV]{Connes-ncg} for details.
Let~$C_b(\mathbb R_+)$ be the space of bounded continuous functions
on $\mathbb R_+$ and $C_0(\mathbb R_+)$ the subspace of
functions vanishing at $\infty$.
Let $\omega$ be  a  positive linear functional
on the quotient space $C_b(\mathbb R_+)/C_0(\mathbb R_+)$ such that
$\omega(1) =1$. For~a positive compact operator $T\in \mathcal L^{1, \infty}$
with eigenvalues $\{\mu_n\}$, extend $\mu_n$ to a step function
on $\mathbb R_+$ and let $M_T(\lam)$ be its  Ces\'aro mean,
which is a bounded continuous function on~$\mathbb R^+$.
The Dixmier trace of $T$ is then defined
by
$$
\tr_{\omega} T=\omega (M_T).
$$
It is then extended to all of $\mathcal L^{1, \infty}$
be linearity. In particular it is bounded and vanishes
on trace class operators. The fact that we will need is that
$$
\tr_{\omega} T=\lim_{N\to \infty} \frac 1{\log N}\sum_{n=1}^{N}\mu_n(T)
$$
if $T$ is a positive operator and if the right hand side exists.

\begin{lemm+}\label{delta-n} For any $c\ge 0$
the operator $(c - \rho(\Delta))^{-d}=\rho(c\delta_0-\Delta)^{-d}$ 
is in the Macaev class $\mathcal L^{1,\infty}$.
\end{lemm+}
\begin{proof}It follows from  (\ref{eq:ossci}) that
the eigenvalues of $(c -\rho(\Delta))^d$ are
$(c+\pi(m +\frac d2))^d$, $m=0, 1, \cdots$, each of  multiplicity
$d_m:=\dim \{w^{\a}, |\alpha|=m\} =\binom{d+m-1}{d-1}\approx m^{d-1}$.
The partial sums thus satisfy
$$
\sum_{m\le N}(c+\pi(m +\frac d2))^{-d} d_m
\approx \sum_{m\le N}(c+\pi(m +\frac d2))^{-d} m^{d-1}
\approx \log N,
$$
completing the proof.
\end{proof} 

\begin{prop+} 
Let $F\in \mathcal {PT}(-2d, 1)$. Then the Weyl transform
$\rho(F)
$ is in the Macaev class $\mathcal L^{1,\infty}$.
\end{prop+}

\begin{proof}
%
%
By~(3.5.6) in~\cite{Howe-jfa80}, 
\begin{equation}
 \hat\Delta=-\frac{\pi^2}4|w|^2,  
 \label{eq:hat-delta}
\end{equation}
so~$-\Delta\in\mathcal{PT}(2,1)$, whence by~(\ref{eq:pridana}) 
$(-\Delta)^{\ast d}\in\mathcal{PT}(2d,1)$ and $(-\Delta)^{\ast d}\ast F
\in\mathcal{PT}(0,1)$. By~the Calder\'on-Vaillancourt theorem
\cite[Theorem~3.1.3]{Howe-jfa80}, the~corresponding Weyl transform,
$\rho(-\Delta)^d\rho(F)$, is~bounded. Hence by the previous lemma 
$\rho(F)\in\mathcal L^{1,\infty}$, since the Macaev class 
$\mathcal L^{1,\infty}$ is an ideal.
\end{proof}

\section{Dixmier trace formula for Toeplitz operators}
\begin{theo+} 
Let $F\in \mathcal {PT}_{rad}(-2d, 1)$ with the
principal symbol $\sig_{-2d}(\hat F)$ as defined in 
(\ref{eq:pr-sym}). Then the Dixmier
trace $\tro\rho(F)$ is independent of $\omega$
and is given by
$$
\tro\rho(F)= \frac{\pi^d}{4^{d}}
\int_{S} \hat \sig_{-2d}(F)(w)
$$
where $\int_S$ is the normalized integral over the unit sphere.
\end{theo+}

\begin{proof} The proof is quite similar to that of Connes \cite{Connes-cmp88}
for pseudo-differential operators on compact manifolds.
Namely, by \cite[Theorem~4.2.5]{Howe-jfa80} and the definition
of~$\mathcal{PT}_{rad}$, the Dixmier trace $\tro\rho(F)$ depends only on 
the leading symbol of $\sig_{-2d}(\hat F)$ and defines a positive measure
on the unit sphere $S$ in $\mathbb C^d$. By the unitary invariance
of $\rho(F)$ the measure has to be a constant multiple of the area measure. 
To~find the constant we 
note that the symbol of $c\delta_0-\Delta$, $c>0$, is absolutely elliptic
in the sense of (4.2.20) in~\cite{Howe-jfa80}, and thus by pp.~246--247
in~\cite{Howe-jfa80} we can construct $F_0\in\mathcal{PT}_{rad}(-2d,1)$ 
such that $\rho(F_0)=(c-\rho(\Delta))^{-d}$. The eigenvalue of $\rho(F_0)$  
on the space of all $m$-homogeneous polynomials~is, by the proof of Lemma~3.1,
$$
\frac 1{(c+\pi(m+\frac d2)  )^d} .
$$
Its~Dixmier trace exists and is
$$
\tro\rho(F_0)=\frac 1{\pi^d}.
$$
On the other hand
the principal symbol $\sig_{-2d}( F_0)$ is the constant function $(4/\pi^2)^d
|w|^{-2d}$ by the definition (cf.~(\ref{eq:hat-delta})), whose integration
over the sphere is $(4/\pi^2)^d$. This completes the proof.
\end{proof}

To apply our result to Toeplitz operators we need to introduce
some more notation.
We let
$$
\partial_j^b =\partial_j - \bar z_j R, \quad
\bar \partial_j^b =\bar \partial_j -  z_j \bar R, 
$$
be the boundary Cauchy-Riemann operators \cite{Rudin-ball},
where 
$R=\sum_{j=1}^d z_j\partial_j$ is the holomorphic radial derivative.
As vector fields they are linearly dependent, to wit,
\begin{equation}
  \label{eq:lin-dep}
 \sum_{j=1}^d z_j\partial_j^b =0, 
\quad  \sum_{j=1}^d \bar z_j\bar \partial_j^b =0. 
\end{equation}

\begin{defi+}
We define a bracket $\{f, g\}_b$ for smooth functions $f$ and $g$ on $S$ by
 $$
\{f, g\}_b := \sum_{j=1}^d (\partial^b_j f \, \dbar^b_j g  -
\dbar^b_j f \, \partial^b_j g)
$$
and call it the boundary Poisson bracket.
\end{defi+}

\begin{lemm+}
 Let $F$ and $G$  be two functions in 
$\mathcal {PT}_{rad}(0, \mu)$
with principal symbols
$$\sig_0(F)(z) = f\Big(\frac z{|z|}\Big),
\quad
\sig_0(G)(z)= g\Big(\frac z{|z|}\Big)
$$
for $f$ and $g$ in $C^\infty(S)$.
Then the principal symbol of $F\ast G- G\ast F$ is given
by 
$$
\sig_{-2}( F\ast G- G\ast F)(z)=
\frac4\pi  \{f, g\}_b (\frac z{|z|})   |z|^{-2} .
$$
\end{lemm+}
\begin{proof} By the general result for the symbol calculus
for pseudo-Toeplitz operators, cf.~(2.2.5) in~\cite{Howe-jfa80},
we~have $F\ast G-G\ast F\in\mathcal{PT}_{rad}(-2\mu,\mu)$ with
the principal symbol 
$$
\sig_{-2}( F\ast G- G\ast F)(z)
=\frac4\pi  \{\sig_0(F)),  \sig_0(G) \}(z) ,
$$
where $\{\cdot,\cdot\}$ is the ordinary Poisson bracket in complex coordinates
$$
\{\Psi, \Phi\}
:= \sum_{j=1}^d (\partial_j \Psi \, \dbar_j \Phi 
- \partial_j \Phi \, \dbar_j \Psi ).
$$
The function $\sig_{-2}( F\ast G- G\ast F)(z)$
is positive homogeneous degree of $-2$. 
We~need only to compute it for $z\in S$. 
Defining 
the Reeb vector field $E$ and the outward normal
vector field $N$ in terms of the
radial derivative $R$,
$$
E:=\frac12(\overline R-R), 
\quad N:=\overline R+R,  $$
we can write
$$ R=-E+\frac N2.
$$
Note
that $E$ and $N$ are all  well-defined on~$S$.
 The~vector field
$\partial^b_j-\overline z_j E$ is thus a well-defined vector field on~$S$, and
for any function $\Phi(z)=\phi(\frac{z}{|z|})$ we have
$$ \partial_j\Phi(z) = (\partial^b_j+\overline z_j R)\Phi(z) 
=(\partial^b_j-\overline z_j E+\frac{\overline z_j}2N)\Phi(z)
=(\partial^b_j-\overline z_j E)\phi(z),  $$
since $N\Phi(z)=0$ by homogeneity. Similarly $\dbar_j\Phi=(\partial^b_j+z_j E)
\phi$ on~$S$.
From this it follows that for $z\in S$
\begin{equation*}
  \begin{split}
\{\sig_0 (F), \sig_0 (G)\}(z)
&= \sum_{j=1}^d \Big((\partial^b_j f(z) - \overline z_j E f (z))
(\dbar^b_j g(z) +  z_j E g(z)) \\
& \qquad\qquad -(\dbar^b_j f(z) - \overline z_j E f (z) )
(\partial^b_j g(z) + z_j E g(z))\Big) \\
&=\{f, g\}_b,    
  \end{split}
\end{equation*}
by using   (\ref{eq:lin-dep}).
\end{proof}

\begin{theo+} Let $f_1, g_1, \cdots, f_d, g_d$ be smooth functions on $S$,
 $\tilde f_1, \tilde g_1, \cdots, \tilde f_d, \tilde g_d$ their 
smooth extensions to $B$ and $T_{\tilde f_1},T_{\tilde g_1},\cdots,
T_{\tilde f_d},T_{\tilde g_d}$ the associated Toeplitz operators on
$\mathcal H_{\nu}$ for $\nu \ge d$. Then the
product $\prod_{j=1}^d [T_{\tilde f_j}, T_{\tilde g_j}]$
is in the Macaev class and its Dixmier trace is given~by
$$
\tro\prod_{j=1}^d [T_{\tilde f_j}, T_{\tilde g_j}]
= \int_S \prod_{j=1}^d \{f_j, g_j\}_b  .
$$ 
\end{theo+}

\begin{proof} The proof is straightforward from the preceding lemma,
formula (2.2.5) in~\cite{Howe-jfa80} and Theorem~4.1.   \end{proof}

We apply our result to Hankel operators with anti-holomorphic
symbols. Let $f$ be a holomorphic function in a neighborhood
of $B$ and $H_{\bar f}g=(I-P)\bar f g$, $g\in \mathcal H_\nu$
the Hankel operator. Then
$$
[T_{\bar f},  T_f]=
[T_f^\ast,  T_f] = |H_{\bar f}|^2 =H_{\bar f}^\ast H_{\bar f}.
$$

\begin{coro+}  Let $f$ be as above. Then the Hankel operator
is in $\mathcal L^{2d, \infty}$, equivalently the commutator 
$[T_{\bar f},  T_f]$ is in $\mathcal L^{d, \infty}$ and we have
$$
\tro |H_{\bar f}|^{2d}
= \tro ([T_{\bar f}, T_f]^{d})
= \int_{S}(|\nabla f|^2-|Rf|^2)^d .
$$
\end{coro+}

Notice that  $H_{\bar f}$ is in the Schatten class $\mathcal L^{p}$
for $p> 2d$ and that its Schatten norm is
$$
 \Vert H_{\bar f}\Vert_{p}^p\approx \int_{B} (1-|z|^2)^p 
(|\nabla f|^2-|Rf|^2)^{\frac p2}dm(z);
$$
see  \cite{AFJP-jlms}
and \cite{Zhu-ajm91} for the Bergman space case ($\nu=d+1$)
and Hardy space ($\nu =d$).  Our result formula provides thus
a limiting result of the above estimates, and it is interesting
to note that estimate has  an equality as its limit for $p\to 2d$.

\def\cprime{$'$} \newcommand{\noopsort}[1]{} \newcommand{\printfirst}[2]{#1}
  \newcommand{\singleletter}[1]{#1} \newcommand{\switchargs}[2]{#2#1}
  \def\cprime{$'$} \def\cprime{$'$} \def\cprime{$'$}
\providecommand{\bysame}{\leavevmode\hbox to3em{\hrulefill}\thinspace}
\providecommand{\MR}{\relax\ifhmode\unskip\space\fi MR }
\providecommand{\MRhref}[2]{%
  \href{http://www.ams.org/mathscinet-getitem?mr=#1}{#2}
}
\providecommand{\href}[2]{#2}

\end{document}